\documentclass[12pt,reqno]{amsart}
\usepackage{fullpage}
\usepackage{mathpple}

\newtheorem{theorem}{Theorem}
\newtheorem{fact}[theorem]{Fact}

\renewcommand{\mod}[1]{{\ifmmode\text{\rm\ (mod~$#1$)}\else\discretionary{}{}{\hbox{ }}\rm(mod~$#1$)\fi}}

\newsymbol\dnd 232D

\newcommand{\Zn}{{\mathbb Z}_n}
\newcommand{\Zp}{{\mathbb Z}_p}

\begin{document}

\title{A simple polynomial for a simple transposition}
\author{Greg Martin} 
\address{Department of Mathematics \\ University of British Columbia \\ Room 121, 1984 Mathematics Road \\ Canada V6T 1Z2}
\email{gerg@math.ubc.ca}
\maketitle

We are used to thinking of polynomials as very special functions, and with good reason when the domain and range are the real numbers or the rational numbers. However, the situation can be counterintuitive over other rings: over finite fields, for example, {\em every} function is a polynomial! It can thus be interesting to look at familiar functions over less familiar rings and find out what sort of polynomial represents them.

Given any ring $R$, we say that the polynomial $P(x)\in R[x]$ {\em represents} the function $f(x)\colon R\to R$ if $P(a)=f(a)$ for every $a\in R$. In a note in this Monthly, Chen and Mullen \cite{CM} considered when various permutations of the elements of $\Zn$ can be represented by polynomials over $\Zn$. One of the remarks they made is that for odd primes $p$, the polynomial
\begin{equation}
f(x) = -\big[ \big( (x-1)^{p-2} + 1 \big)^{p-2} -1\big]^{p-2}
\label{CM.poly}
\end{equation}
represents the transposition $(0\; 1)$ over $\Zp$, that is,
\begin{equation}
f(0)\equiv1\mod p, \quad f(1)\equiv0\mod p, \quad f(a)\equiv a\mod p \text{ for all }2\le a\le p-1.
\label{trans.conditions}
\end{equation}
(For the even prime $p=2$ we can take simply $f(x)=1-x$.)
Although this is certainly true, we know (see Fact \ref{zpzp} below) that any function from $\Zp$ to $\Zp$ can be represented uniquely by a polynomial of degree at most $p-1$. Since the polynomial \eqref{CM.poly} has degree $(p-2)^3$, it is not in this canonical form when $p\ge5$.

As it happens, the natural question of which polynomial of degree at most $p-1$ represents the transposition $(0\; 1)$ over $\Zp$ has quite a nice answer:
\begin{theorem}
For all odd primes $p$, the polynomial
\begin{equation}
f(x) = x^{p-2} + x^{p-3} + \dots + x^3 + x^2 + 2x + 1  \label{my.poly}
\end{equation}
represents the transposition $(0\; 1)$ over $\Zp$, that is, the congruences \eqref{trans.conditions} are all satisfied.
\label{main.thm}
\end{theorem}

By a linear change of variables, we see that any transposition $(a\; b)$ of elements of $\Zp$ is represented by an analogous polynomial in $\Zp[x]$, namely
\begin{equation}
(b-a) \bigg( \bigg( \frac{x-a}{b-a} \bigg)^{p-2} + \bigg( \frac{x-a}{b-a} \bigg)^{p-3} + \dots + \bigg( \frac{x-a}{b-a} \bigg)^3 + \bigg( \frac{x-a}{b-a} \bigg)^2 + 2\bigg( \frac{x-a}{b-a} \bigg) + 1 \bigg) + a.
\label{more.general}
\end{equation}

To prove Theorem \ref{main.thm}, we need the well-known fact alluded to earlier, for which we give two proofs, each standard.

\begin{fact}
Every function $f:\Zp\to\Zp$ is represented by a unique polynomial of degree at most $p-1$.
\label{zpzp}
\end{fact}

\begin{proof}
To prove the claim of uniqueness, we notice that two polynomials $f_1(x)$ and $f_2(x)$ represent the same function over $\Zp$ if and only if $(f_1-f_2)(a) \equiv0\mod p$ for every integer $a$. That means the polynomial $(f_1-f_2)(x)\in\Zp[x]$ is divisible by $x-a$ for every $0\le a\le p-1$. Since $\Zp$ is a field, $\Zp[x]$ is a unique factorization domain, and so $f_1-f_2$ must be divisible by the product $\prod_{a=0}^{p-1} (x-a)$. In particular, this divisor has degree $p$, so if $f_1$ and $f_2$ have degree at most $p-1$, then they must be equal.

It is therefore enough to prove every function from $\Zp$ to $\Zp$ is represented by at least one polynomial of degree at most $p-1$.

\smallskip{\em Proof 1:} Since $t^{p-1}\equiv0\text{ or }1\mod p$ according to whether or not $t\equiv0\mod p$ by Fermat's little theorem, we can write
\begin{equation}
f(x) \equiv \sum_{a=0}^{p-1} f(a) \big( 1-(x-a)^{p-1} \big)\; \mod p,
\label{each.coeff}
\end{equation}
and the right-hand side is clearly a polynomial of degree at most $p-1$.

\smallskip{\em Proof 2:} There are $p^p$ polynomials of degree at most $p-1$, and by the uniqueness proved above, each one represents a different function. But there are exactly $p^p$ functions from $\Zp$ to $\Zp$ as well. Therefore every function must be represented by such a polynomial.
\end{proof}

Although we didn't need to know it in the proof of Fact \ref{zpzp}, it happens that the product $\prod_{a=0}^{p-1} (x-a)$ is equal (in $\Zp[x]$) to the simple polynomial $x^p - x$. This is true because every integer $a$ is a root of the latter polynomial by Fermat's little theorem, and by unique factorization we again conclude that $\prod_{a=0}^{p-1} (x-a)$ divides $x^p - x$ in $\Zp[x]$; since the two polynomials have the same degree and leading coefficient, they must be equal. By dividing both sides by $x$, it follows also that
\begin{equation}
\prod_{a=1}^{p-1} (x-a) = x^{p-1}-1  \label{fltl}
\end{equation}
as polynomials in $\Zp[x]$.

In fact, we have shown that two polynomials in $\Zp[x]$ represent the same function if and only if they differ by a multiple of $x^p-x$. Certainly it isn't obvious that the difference of the two polynomials \eqref{CM.poly} and \eqref{my.poly} has $x^p-x$ as a factor! Of course, this conclusion relies on Theorem \ref{main.thm}, and so we prove it now.

\begin{proof}[Proof of Theorem \ref{main.thm}]
Notice that
\[
x^{p-2} + x^{p-3} + \dots + x + 1 = \frac{x^{p-1}-1}{x-1} = \prod_{a=2}^{p-1} (x-a)
\]
(the first equality holds over any ring $R$, while the second, a rearrangement of \eqref{fltl}, is particular to $\Zp$). Therefore the left-hand side yields 0\mod p when $x$ equals any of the integers $a=2,3,\dots,p-1$. If we now define
\[
f(x) = (x^{p-2} + x^{p-3} + \dots + x + 1) + x = x^{p-2} + x^{p-3} + \dots + 2x + 1,
\]
then clearly $f(a)\equiv 0+a=a\mod p$ for all $2\le a\le p-1$. However, the values $f(0) = 1$ and $f(1) = p \equiv 0\mod p$ are easy to calculate. We see that $f(x)$ represents the transposition $(0\; 1)$ as claimed.
\end{proof}

We remark in passing that the polynomial \eqref{more.general}, which represents the transposition $(a\; b)$ modulo $p$, can also be written in the form
\begin{equation}
(b-a)^2 \frac{x^p-x}{(x-a)(x-b)} + x.  \label{other.form}
\end{equation}
It is not hard to show algebraically that the two forms are equivalent, particularly if we write $x^p-x = (x-a)^p - (b-a)^{p-1}(x-a)$ using two applications of Fermat's little theorem. However, the expression \eqref{other.form} can be seen to represent the transposition $(a\; b)$ immediately, by using a finite-field analogue of l'H\^opital's rule: if $f(c)\equiv 0\mod p$ then
\[
\frac{f(x)}{x-c}\bigg|_{x=c} \equiv f'(c)\mod p,
\]
which follows from expanding the polynomial $f$ as a Taylor series about $x=c$.

Since functions on $\Zp$ are canonically represented by polynomials of degree at most $p-1$, the polynomials that represent transpositions, having strictly lower degree, are somewhat special. In hindsight we could have predicted that these polynomials would have degrees less than $p-1$. The following fact is rather less well-known than Fact \ref{zpzp}; again we give two proofs, but this time the second proof seems to be new.

\begin{fact}
The polynomial of degree at most $p-1$ that represents the function $f(x)$ has degree at most $p-2$ if and only if
\begin{equation}
\sum_{a=0}^{p-1} f(a) \equiv 0\mod p.
\label{topcong}
\end{equation}
\end{fact}

\smallskip{\em Proof 1:} Expanding the right-hand side of \eqref{each.coeff}, we see that the coefficient of $x^{p-1}$ is exactly $\sum_{a=0}^{p-1} f(a)$.

\smallskip{\em Proof 2:} There are exactly $p^{p-1}$ functions from $\Zp$ to $\Zp$ satisfying the congruence \eqref{topcong} and exactly $p^{p-1}$ polynomials in $\Zp[x]$ of degree at most $p-2$. It thus suffices, as before, to prove that every function $f$ satisfying \eqref{topcong} is represented by some polynomial of degree at most $p-2$. Given such a function $f$, for every positive integer $n$ define
\[
F(n) = \sum_{a=1}^n f(a).
\]
Because of the assumption \eqref{topcong} on $f$, this function $F$ is actually a function from $\Zp$ to $\Zp$: to see that it is periodic with period $p$, we calculate
\[
F(n+p)-F(n) = \sum_{a=n+1}^{n+p} f(a) \equiv \sum_{a=0}^{p-1} f(a) \equiv 0\mod p,
\]
since the values $f(n+1),\dots,f(n+p)$ are a rearrangement of $f(0),\dots,f(p-1)$ modulo $p$. Therefore $F$ is represented by a polynomial $P(x)$ of degree at most $p-1$ by Fact \ref{zpzp}. But since $f(x) = F(x)-F(x-1)$, we see that $f$ is represented by the polynomial $P(x)-P(x-1)$, which has degree one less than the degree of $P$ itself, hence at most $p-2$.\qed\smallskip

Both proofs generalize (the first one more readily) to show that the polynomial representing the function $f$ has degree equal to $p-1-k$ if and only if
\begin{multline*}
\sum_{a=0}^{p-1} f(a) \equiv0\mod p,\quad \sum_{a=0}^{p-1} af(a) \equiv0\mod p,\quad \dots, \\
\sum_{a=0}^{p-1} a^{k-1}f(a) \equiv0\mod p,\quad \sum_{a=0}^{p-1} a^kf(a) \not\equiv0\mod p.
\end{multline*}

Note that the polynomial representing any permutation $\sigma$ of the elements of $\Zp$ (not only a transposition) has degree at most $p-2$, since $\sum_{a=0}^{p-1} \sigma(a) \equiv \sum_{a=0}^{p-1} a\equiv0\mod p$. The possible degrees of polynomials representing permutations have been studied in detail; it can be shown, for example, that the degree of such a polynomial can never divide $p-1$. The interested reader can refer to Lidl and Niederreiter \cite[Chapter 7]{LN} to learn more about this topic.

\bigskip{\small
{\it Acknowledgements.} L'auteur remercie le Centre de recherches math\'ematiques \`a l'Universit\'e de Montr\'eal de son hospitalit\'e pendant la composition de cet article. The author was supported in part by a Natural Sciences and Engineering Research Council grant.}

\bibliographystyle{amsplain}
\bibliography{permpoly}

\end{document}